\documentclass[twoside,leqno]{article}

\usepackage[letterpaper]{geometry}

			\usepackage{xcolor}
			
			\usepackage[urlcolor=red,colorlinks=true]{hyperref}

\usepackage{amssymb}

\usepackage{siamproceedings}

\usepackage[T1]{fontenc}
\usepackage{amsfonts}
\usepackage{graphicx}
\usepackage{epstopdf}
\usepackage{enumitem}
\usepackage{algorithmic}
\ifpdf
  \DeclareGraphicsExtensions{.eps,.pdf,.png,.jpg}
\else
  \DeclareGraphicsExtensions{.eps}
\fi

\newsiamremark{remark}{Remark}
\newsiamremark{hypothesis}{Hypothesis}
\crefname{hypothesis}{Hypothesis}{Hypotheses}
\newsiamthm{claim}{Claim}

\usepackage{amsopn}

\providecommand{\keywords}[1]
{
  \small	
  \textbf{\textit{Keywords---}} #1
}

\newtheorem{conjecture}[theorem]{Conjecture}

\newcommand{\bprop}{\begin{proposition}}
\newcommand{\ele}{\end{lemma}}
\newcommand{\ecor}{\end{corr}}
\newcommand{\edeff}{\end{deff}}

\newcommand{\eprop}{\end{proposition}}

\newcommand{\Rn}{{\mathbb R}^n}

\newcommand{\la}{\lambda}

\newcommand{\e}{\varepsilon}

\renewcommand{\Pi}{\varPi}

\renewcommand{\epsilon}{\varepsilon}

\newcommand{\Rt}{{\Bbb R}^3}

\newcommand{\parital}{\partial}

\newcommand{\R}{{\mathbb R}}
\newcommand{\Td}{{\mathbb T}_\delta}
\newcommand{\Tr}{{\mathbb T}_\rho}

\newcommand{\mC}{{\mathcal C}}

\newcommand{\one}{{\bf 1}}

\begin{document}

\title{\Large The Work of Hong Wang}
    \author{Christopher D. Sogge\thanks{Department of Mathematics, Johns Hopkins University (\email{sogge@jhu.edu})  The  author was  supported in part by the NSF (DMS-2348996.}}
    
    \date{}
    
   %\thanks{The  author was  supported in part by the NSF (DMS-2348996).}

%\subjclass[2020]{Primary 42B25, 42B20; Secondary 35L05}

\maketitle

\begin{abstract} Hong Wang has made several tremendous contributions to harmonic
 analysis and geometric
measure theory.  On the 23rd of July 2026, she was awarded the Fields Medal.  In this short article
we shall focus on her work on the Kakeya problem and other 
problems related to the Fourier transform.  We will try to put
her work in historical context and present a small selection of her contributions to  harmonic analysis.
\end{abstract}

\keywords
{Fourier transform, Kakeya sets, local smoothing, multiscale analysis}

%\subjclass[2020]{Primary 42B25, 42B20; Secondary 35L05}

% Copyright Statement
% When submitting your final paper to a SIAM proceedings, it is requested that you include
% the appropriate copyright in the footer of the paper.  The copyright added should be
% consistent with the copyright selected on the copyright form submitted with the paper.
% Please note that "20XX" should be changed to the year of the meeting.

% Default Copyright Statement

%\fancyfoot[R]{\scriptsize{Copyright \textcopyright\ 20XX by SIAM\\
%Unauthorized reproduction of this article is prohibited}}

% Depending on which copyright you agree to when you sign the copyright form, the copyright
% can be changed to one of the following after commenting out the default copyright statement
% above.

%\fancyfoot[R]{\scriptsize{Copyright \textcopyright\ 20XX\\
%Copyright for this paper is retained by authors}}

%\fancyfoot[R]{\scriptsize{Copyright \textcopyright\ 20XX\\
%Copyright retained by principal author's organization}}

%\pagenumbering{arabic}
%\setcounter{page}{1}%Leave this line commented out.

 		\cfoot{\thepage}

\section{Introduction.}
Hong Wang received the Fields Medal for her stellar work in harmonic analysis and geometric measure theory.  She  resolved the
three-dimensional Kakeya conjecture, the planar local smoothing conjecture for the wave equation and the Furstenberg set conjecture.
Each of these conjectures was several decades old and had attracted the attention of the top researchers in the field.  She has also
been able to obtain the best known results for several other longstanding problems, including Stein's restriction problem for the Fourier
transform, the Bochner-Riesz problem concerning the convergence of truncated Fourier transforms, as well as the Falconer distance
problem relating Hausdorff dimension and the structure of distance distributions.

A major theme of Wang's work is the novel use of multiscale methods in combination with deep and innovative geometric combinatorial
arguments.  Often these involve the discovery and exploitation of new geometric structures.  Wang has also had remarkable success
in extending and refining modern techniques in harmonic analysis and geometric measure theory, especially decoupling theory and
multiscale analysis.  The methods that she has introduced in her breakthroughs on several central problems in analysis have already
been impactful in related areas such as geometry and partial differential equations.

Much of Wang's work, including her spectacular work with Joshua Zahl,  \cite{WZ2}, \cite{WZ3}, \cite{WZ1}, resolving the three-dimensional
Kakeya conjecture and the remarkable work with Larry Guth and Ruixiang Zhang \cite{GWZls} settling the two-dimensional local smoothing
conjecture for the wave equation, is related to the Euclidean Fourier transform
\begin{equation}\label{1.1}
\Hat f(\xi)= \int _{{\mathbb R}^n} e^{-2\pi i x\cdot \xi} \, f(x) \, dx.
\end{equation}
Indeed, the Kakeya conjecture and the local smoothing conjecture belong to  a family of conjectures, a ``tower'' of conjectures
(see \cite{Tow}), each of which has deep and surprising connections with the Fourier transform and one another.
In this article we shall focus on Wang's work on these problems.  Before stating her breakthroughs, we shall give some historical
background which helps to emphasize the significance of her results.

\section{Background: The Kakeya Problem}

Before the results of Wang and Zahl  \cite{WZ2}, \cite{WZ3}, \cite{WZ1} there were two other milestones that we shall 
describe in this section.

\subsection{Early Kakeya Results}

For many years people have been interested in the ``size'' of Kakeya sets, which we define as follows:

\begin{definition}
A compact subset ${\mathcal K}\subset \Rn$, $n\ge 2$, is called a {\em Kakeya set} if it contains a unit-line segment pointing
in every direction.
\end{definition}

This assumption means that, if $S^{n-1}\subset \Rn$ is the unit sphere, then for each $\omega \in S^{n-1}$ one  can find
$x_\omega\in \Rn$ so that $\gamma_\omega =\{ \, x_\omega +t\omega: \, \, t\in [0,1]\}\subset {\mathcal K}$.

Without knowing of each other’s work, in 1917, Abram Besicovitch and S\=oichi Kakeya were the first to study Kakeya sets.
This was due to the relative isolation of Russian mathematicians at that time on account of the Russian revolution.

They had very different motivations for studying Kakeya sets.  Besicovitch was working on a problem in Riemannian integration
concerning the extent to which the Fubini theorem for the Lebesgue integral extends to the Riemannian setting.  Specifically,
as stated in Besicovitch's survey paper \cite{Be63}:

{\em Given a function of two variables, Riemann-integrable on a plane domain, does there always exist a pair of mutually perpendicular
directions such that the repeated simple integration along the two directions exists and gives the value of the integral over
the domain?}

Besicovitch realized\footnote{For the reasons we have stated before, unfortunately, Besicovitch's 1917 work is lost to time and we are
unable to provide a reference.}
 that if there were planar Kakeya sets of Lebesgue measure zero then these would provide a counterexample, 
ensuring that the answer must be ``no''.  See \cite{Fal} for details.

As we shall see, this is a recurrent theme.  ``Small'' Kakeya sets yield negative results in analysis, and, consequently, one might hope that
Kakeya sets cannot be ``too small.''

Kakeya's motivation was much different.  In his 1917 paper, ``Some problems on maximum and minimum regarding ovals,'' he proposed
54 problems.  The most famous is what we now call the Kakeya needle problem:

{\em
In the class of figures in which a segment of length $1$ can be turned around through $360^\circ,$ remaining always within the figure,
which one has the smallest area?
}

\begin{figure}[htbp]
  \centering
  \includegraphics[width=0.6\columnwidth]{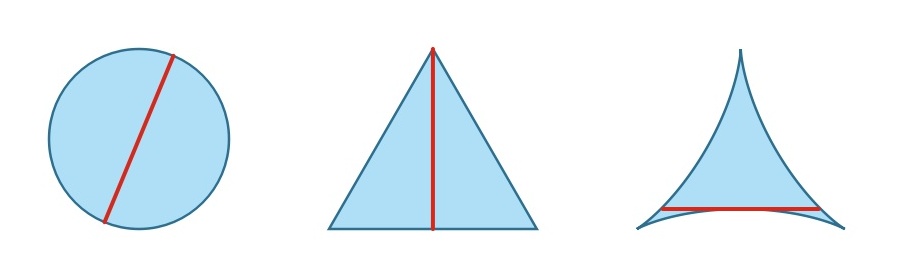}
  \caption{Revolvable sets.}
  \label{fig1}
\end{figure}

Kakeya called sets with this property ``revolvable.''  Thus, his needle problem asks one to find a planar revolvable set of minimal area.
He gave three examples: discs of diameter one, equilateral triangles of height one and deltoids\footnote{Recall that a deltoid is a roulette curve.  One of unit height is formed by marking a circle of diameter $1/2$ and then rolling it within
a circle of diameter $3/2$ to trace out the deltoid via the marked point.}
of height one, as shown in Figure~\ref{fig1}.

The three sets in the figure have decreasing area.  The disc of course is of area
$\pi/4 \approx .78$, the triangle $\sqrt{3}/3\approx .58$, while the deltoid has area
$\pi/8\approx .39$.  Thus, the deltoid has half the area of the disc.  Furthermore,
when executing a three-point turn one's car basically sweeps out a deltoid.

Perhaps motivated by these sorts of considerations Kakeya conjectured that the deltoid
is the solution of his needle problem, i.e., a revolvable set of minimal area.  Unaware of 
Besicovitch's work \cite{Be19}, the Kakeya needle problem aroused much interest.  For instance, 
in 1925, Birkhoff  writing about unsolved problems in his book, ``The Origin, Nature and Influence of Relativity'' \cite{Bir}, first mentions
the four-color problem and then added: {\em ``Of like intriguing simplicity is the question raised a few years ago
by the Japanese mathematician Kakeya.''}

Besicovitch established the following milestone which provided a negative solution to his problem concerning the Riemann integral,
as well as the Kakeya needle problem: 
\begin{theorem}[Besicovitch 1919/1928]\label{thmBe}
There are Kakeya sets  ${\mathcal K}\subset \Rn$ of {\em measure zero}, $|{\mathcal K}|=0$.  
Additionally, there are planar 
revolvable
sets of {\em arbitrarily small area}.
% in which a unit line
%segment can be continuously rotated through 180$^\circ$ while remaining in the set.
\end{theorem}

Besicovitch's 1919 paper \cite{Be19} establishing zero-measure planar Kakeya sets appeared in
an obscure Russian journal that went defunct after only publishing a couple of  volumes due to the turbulent times
in Russia.  At the time, Besicovitch was unaware of Kakeya's problem, and mathematicians outside
of Russia were largely unaware of Besicovitch's work.  Consequently, he was led to republish his construction
of zero-measure Kakeya sets (often called ``Besicovitch sets'') in his 1928 paper
\cite{Be28}.  In this paper he added a one-paragraph modification of his construction,
due to the Hungarian mathematician P\'al\footnote{P\'al also showed that the equilateral triangle of
height one is the smallest convex revolvable set in \cite{Pal}.},
showing that there are revolvable sets of arbitrarily small area.  Additionally, Besicovitch was only
interested in planar sets, but one can obtain zero-measure Kakeya sets in $\Rn$, $n\ge 3$, 
by taking the Cartesian product of a zero-measure planar set ${\mathcal K}$ and
$[0,1]^{n-2}$.

\begin{figure}[htbp]
  \centering
  \includegraphics[width=1.0\columnwidth]{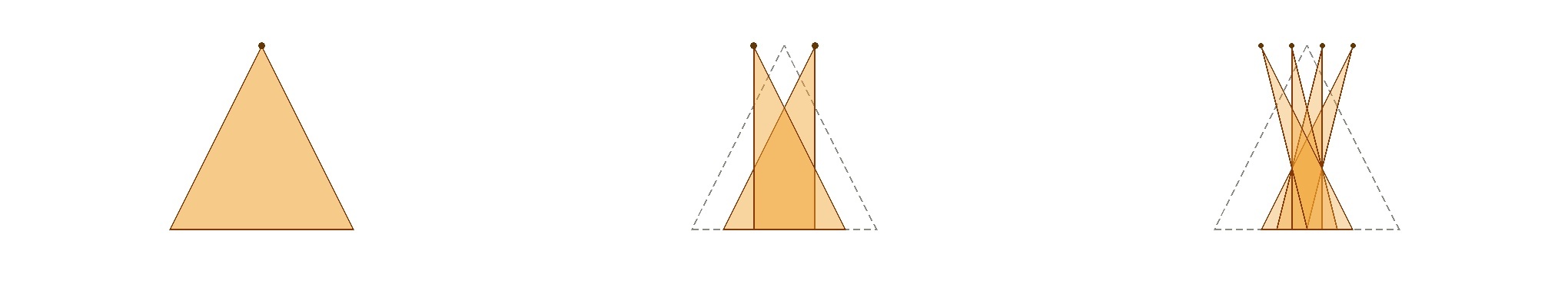}
  \caption{Besicovitch's construction: ${\mathcal K}_0, \, \, {\mathcal K}_1, \, \, {\mathcal K}_2$}
  \label{fig2}
\end{figure}

Besicovitch's construction
is based on the observation that if one bisects the top vertex of a triangle and then
slides the two halves along the bottom side so that they overlap, then the resulting object will have smaller area.  
If one starts with an equilateral triangle as in Figure~\ref{fig2}, the resulting set will include a unit line segment
pointing in each direction swept out by the top vertex.  Thus, since segments have two directions, it contains
segments pointing in $120^\circ = 60^\circ + 60^\circ$ worth of directions.  If one repeats this construction,
as in Figure~\ref{fig2}, the resulting set has even smaller area but still contains a unit line segment pointing
in all of these directions.  
By iterating this procedure, after $\ell$ steps, Besicovitch obtained a set ${\mathcal K}_l$ consisting of $2^\ell$ ``ears'' and having total
area $|{\mathcal K}_\ell|\approx 1/\ell$.
He then set up a limiting argument
to obtain a zero-measure set containing a unit line segment pointing in all the directions swept out by the top
vertex of the equilateral triangle in Figure~\ref{fig2}.  Of course, by suitably rotating two copies of this set
and taking unions with the original set, one obtains the zero-measure Kakeya set ${\mathcal K}$ in Theorem~\ref{thmBe} when $n=2$.

The argument we have sketched is a slight simplification of Besicovitch's original
argument, due to Perron~\cite{Perr} and  Schoenberg~\cite{Schn}.  Additionally, the sets ${\mathcal K}_\ell$
described above are called ``Perron trees''.

For more details about Besicovitch's construction and the early history of the Kakeya problem, we refer
the reader to his excellent survey article \cite{Be63}.  This paper was written to accompany a delightful  
{\color{red}\href{https://tinyurl.com/2huh4y84}{film}}
funded by a 1958 National Science Foundation grant established
to ``explore the possibilities of mathematical motion pictures''.   This  film about Kakeya sets
was narrated by Besicovitch, 
and it
was one of four mathematical films funded by this early NSF grant.

\subsection{Kakeya Sets and Fourier Analysis}

Another old problem concerns the $L^p$-convergence of the truncated
Fourier integrals,
\begin{equation}\label{2.1}
S_Nf(x)=\int_{\{\xi\in \Rn: \, |\xi|\le N\}} 
e^{2\pi i x\cdot \xi} \, \Hat f(\xi) \, d\xi,
\end{equation}
with the Fourier transform $\Hat f$ defined as in \eqref{1.1}.
Recall that Fourier's inversion formula says that
\begin{equation}\label{2.2}
f(x)=\int_{\Rn} e^{2\pi i x\cdot \xi} \, \Hat f(\xi) \, d\xi,
\end{equation}
for well-behaved functions $f$, for instance $f$ in ${\mathcal S}(\Rn)$,
the space of Schwartz-class functions.  Recall also that, from this and
a density argument, one obtains that there is $L^2$-convergence
of the ball multiplier operators, i.e., 
\begin{equation}\label{2.3}
S_Nf\to f \quad \text{in } \, \, L^2(\Rn) \, \, \text{as } \, \, N\to \infty.
\end{equation}
By this we of course mean that if $f\in L^2(\Rn)$, then
$\|S_Nf-f\|_{L^2(\Rn)}\to 0$ as $N\to \infty$.

The $S_N$ are called ball multiplier operators since the Fourier transform of $S_Nf$ is just
$\Hat f$ multiplied by the indicator function of the $N$-ball, $\{\xi\in \Rn: \, |\xi|\le N\}$.  Also, by
a simple dilation argument
\begin{equation}\label{2.4}
\|S_N\|_{p\to p}=\|S\|_{p\to p}, \quad \text{if } \, \, \, S=S_1,
\end{equation}
with $\|\, \cdot \, \|_{p\to p}$ denoting the $L^p(\Rn)\to L^p(\Rn)$ operator norm.

Given \eqref{2.3}, it is natural to ask whether there is also $L^p(\Rn)$-convergence of the 
truncated Fourier integrals $S_Nf$ in \eqref{2.1}.  In other words, given a 
dimension $n=1,2,\dots$,  for which exponents $p\in (1,\infty)$ do we have
\begin{equation}\label{2.5}
S_Nf\to f \quad \text{in } \, \, L^p(\Rn) \, \, \text{as } \, \, N\to \infty \,  ?
\end{equation}
It is not difficult to see that $1<p<\infty$ is a necessary condition.

This classical problem is called the multiplier problem for the ball.  By the above considerations,
there is an affirmative answer for a given $p$ if and only if $S: L^p(\Rn)\to L^p(\Rn)$, with,
as in \eqref{2.4}, $S$ being the multiplier operator associated with the unit ball, that is, $S=S_1$.

In 1928, M. Riesz~\cite{Rie} completely solved the one-dimensional ball multiplier problem, showing
that
\begin{equation}\label{2.6}
S_Nf\to f \quad \text{in } \, \, L^p(\R) \quad \text{for every} \, \, \, 1<p<\infty.
\end{equation}
Riesz's proof of this landmark result introduced important new tools in harmonic analysis, including
interpolation theory.

In the 1970s Fefferman began his study of the ball multiplier problem in higher dimensions $n\ge2$.  
In \cite{Fball} 
he came
up with the following very surprising full resolution:
\begin{equation}\label{2.7}
\text{If } \, \, n\ge 2\, \, \text{then } \, \, S_Nf \to f \, \,
\text{in } \, \, L^p(\Rn) \, \, \text{if and only if } \, \, p=2. 
\end{equation}
He proved this milestone result by using Besicovitch's iterative construction of zero-measure Kakeya sets mentioned before.
In \cite{Fball}, Fefferman called his {\em Ball Multiplier Theorem}, \eqref{2.7}, an ``unfortunate fact.''

Let us give a very rough sketch of Fefferman's seminal proof linking the Fourier transform and Kakeya sets.  To prove
\eqref{2.7} for, say, $n=2$, it suffices, by duality and the above considerations, to show that
\begin{equation}\label{2.8}
S \, \, \, \text{is unbounded on } \, \, L^p(\R^2) \, \, \, \text{for all }  \, \, p>2.
\end{equation}

To prove \eqref{2.8} Fefferman introduced a remarkable construction relying on the properties of the iterates,
${\mathcal K}_\ell$, described earlier in Besicovitch's construction.  To work spatially within sets of unit-scale
near ${\mathcal K}_\ell$, after recalling
 \eqref{2.4}, one deduces that, to show $S$ is unbounded
on $L^p(\R^2)$ for each $p>2$, it suffices to show that
\begin{equation}\label{2.81}
\|S_{N^2}f\|_{L^p(\R^2)}/\|f\|_{L^p(\R^2)}\ge C_N,
\, \, \text{some } \, \, f=f_N, \, \,
\text{with } \, \, C_N\to \infty, \, \, \text{as } \, \, N\to \infty.
\end{equation}

\begin{figure}[htbp]
  \centering
  \includegraphics[width=0.5\columnwidth]{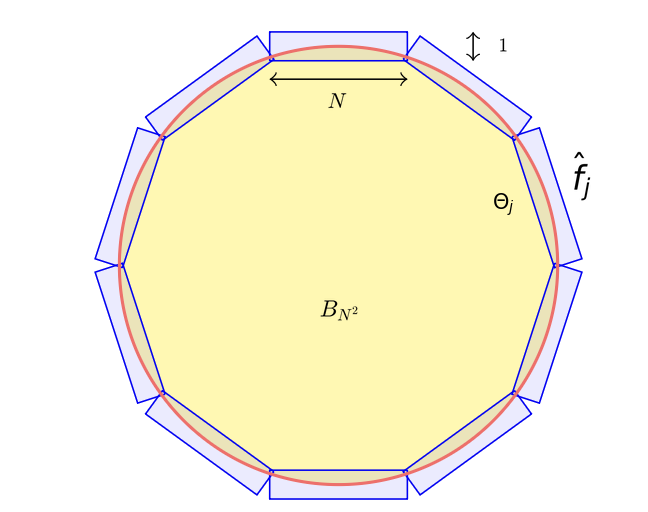}
  \caption{Fefferman's construction.}
  \label{fig3}
\end{figure}

To prove this Fefferman introduces $f=\sum f_j$ with the sum having $N=2^\ell$ terms.  
He chooses each $f_j$ to be the indicator of a $1\times N^{-1}$ rectangle $R_j$ times an
oscillatory (modulation) factor
$$f_j(x)= e^{2\pi i N^2 x\cdot \theta_j} \, \one_{R_j}(x),
$$
with $\theta_j \in S^1$ pointing in the direction of the long side of $R_j$ (the one of length 1).  It then follows
that the Fourier transform $\Hat f_j$ is concentrated on a $1\times N$  rectangle, $\Theta_j$, centered
at $N^2\theta_j$, a point on the boundary of the $N^2$-disc, $B_{N^2}$, with short side pointing
in the direction $\theta_j$ as in Figure~\ref{fig3}.
Note that roughly ``half'' of $\Theta_j$ is inside the $N^2$-disc and roughly ``half'' is outside.
In this construction, only terms with $\theta_j$ corresponding to the direction of the $N=2^\ell$ ``ears'' of ${\mathcal K}_\ell$
are included in the sum defining $f$.

The rectangles $R_j$ chosen will be disjoint, which ensures that
\begin{equation}\label{2.82}
\|f\|_{L^q(\R^2)}\approx 1 \quad \text{all } \, \, q\ge 2.
\end{equation}
The key and difficult step in Fefferman's construction is showing that  $R_j$ can be chosen near each
of the $N=2^\ell$ ``ears'' of ${\mathcal K}_\ell$ so that $S_{N^2}f$ has nontrivial $L^2$-concentration on ${\mathcal K}_\ell$:
\begin{equation}\label{2.83}
\|S_{N^2}f\|_{L^2({\mathcal K}_\ell)} \ge c_0>0, \quad N=2^\ell, \, \ell=1,2,\dots.
\end{equation}
This is a delicate uniform lower bound due to the uncertainty principle and other considerations.
It relies on the fact that the sets ${\mathcal K}_\ell$ in the Besicovitch construction contain rectangles
of size $\sim 1\times N^{-1}$ pointing in each of the $\theta_j$-directions mentioned before, as well as the
fact that $\Hat f_j$ is concentrated on the $1\times N$ rectangle $\Theta_j$ ``half'' of which is inside
$B_{N^2}$ and ``half'' outside.

To use \eqref{2.82} and the difficult Kakeya mass compression estimate \eqref{2.83}, Fefferman just
uses H\"older's inequality and the assumption that $p>2$ to deduce that
$$c_0\le |{\mathcal K}_\ell|^{\frac12-\frac1p}\, \|S_{N^2}f\|_{L^p({\mathcal K}_\ell)}
\le |{\mathcal K}_\ell|^{\frac12-\frac1p}\, \|S_{N^2}f\|_{L^p(\R^2)}.$$
Since Besicovitch's construction yields $|{\mathcal K}_\ell|\approx 1/\ell =1/\log N$, we deduce from this
that
$$\| S_{N^2}f\|_{L^p(\R^2)}/\|f\|_{L^p(\R^2)} \gtrsim  (\log N)^{\frac12-\frac1p}, \quad \text{if } \, \, p>2,
$$
which yields \eqref{2.81}.

As with Besicovitch's theorem, straightforward arguments show that
the higher dimensional $n\ge3$ version of Fefferman's theorem, \eqref{2.7}, follows from the two-dimensional argument.

This is a very rough sketch of Fefferman's incredible observation that there are deep connections between properties
of Kakeya sets and problems involving the Fourier transform, which has ben an important theme in harmonic
analysis since his seminal paper \cite{Fball}.  For more details, the reader can go to Fefferman's self-contained 
paper or Stein's book \cite{SteinH}.
%or the recent survey article \cite{Hick} by Hickman about the Kakeya problem.

Due to  Fefferman's ``unfortunate fact'', \eqref{2.8}, in order to possibly have $L^p$-convergence, $p\ne 2$, of 
truncated Fourier integrals representing $f$, it is necessary to mollify the ball multiplier operators.  A standard
way of doing this is by using the Bochner-Riesz operators $S^\delta_N$, $\delta > 0$, defined by
\begin{equation}\label{2.9}
S^\delta_N f(x) = \int_{\Rn} e^{2\pi ix\cdot \xi} \, \Hat f(\xi) \, \one_{B_N}(\xi) \,  \bigl(1-|\xi/N|^2\bigr)^\delta \, d\xi,
\end{equation}
with $\one_{B_N}$ being the indicator function of the $N$-ball $\{\xi\in \Rn: \, |\xi|\le N\}$.  The analog of the 
ball multiplier operator then becomes, given $p$ and $n$, when do we have
$S^\delta_N f\to f$ in $L^p(\Rn)$?  Since the Fourier multipliers $\one_{B_N}(\xi) \, (1-|\xi/N|^2)^\delta$
become smoother as $\delta$ increases, there is naturally a minimum value of $\delta$ for an affirmative answer.  Due
to \eqref{2.8} and much more straightforward considerations, an optimal
result would be one form of the Bochner-Riesz conjecture
\begin{equation}\label{2.10}
S^\delta_N f\to f \quad \text{in } \, \,
L^p(\Rn) \, \, 
\, \text{if } \, \, \delta>\delta(p,n)=\max \bigl\{ n| \tfrac12 -\tfrac1p |-\tfrac12, \, 0\bigr\}.
\end{equation}
Analogous to the situation for the ball multiplier operators, one has \eqref{2.10} for a given exponent $p$  if and only if
$S^\delta: L^p(\Rn)\to L^p(\Rn)$  with
\begin{equation}\label{BR}
S^\delta f(x) = \int_{\{\xi\in \Rn: \, |\xi|\le 1\}} e^{2\pi ix\cdot \xi} \, (1-|\xi|^2)^\delta \, \Hat f(\xi)\, d\xi.
\end{equation}
%Note that $\delta(p,n)=0$ if $p=2n/(n-1)$.  
%Using this observation and an interpolation argument, one finds 
As a result, one sees that, if $\delta(p,n)$ is as above, then
the Bochner-Riesz convergence in \eqref{2.10} is valid if and only if the standard form of the
Bochner-Riesz conjecture holds:
\begin{equation}\label{2.11} 
S^\delta: L^p(\Rn)\to L^p(\Rn) \, \quad 
\text{if } \, \, \, 
\delta>  \delta(p,n).
\end{equation}

Just like Besicovitch, Fefferman proved his ball multiplier theorem by using the fact that there are zero-measure
Kakeya sets ${\mathcal K}\subset \Rn$, $n\ge2$.  For the Bochner-Riesz conjecture to be valid, it is important
that zero-measure Kakeya sets cannot be ``too small.''  Indeed, if $\text{dim}$ denotes either Minkowski or
Hausdorff dimension, Fefferman's argument also shows that, if there are Kakeya sets ${\mathcal K}\subset \Rn$
satisfying $\text{dim }{\mathcal K}\le n-\alpha$ for some $\alpha>0$, then the bounds in \eqref{2.11} 
cannot be valid for exponents $p\in (p'(\alpha,n),p(\alpha,n)) \, \backslash  \, \{2\},$  for some exponent
$p(\alpha,n)>2$,
depending on $n$ and $\alpha$.  Here $p'(\alpha,n)<2$ denotes conjugate exponent
for $p(\alpha,n)$  defined by $1/p(\alpha,n)+1/p'(\alpha,n)=1$.
Thus,  {\em the 
Bochner-Riesz conjecture implies that Kakeya sets must have full dimension}.

Around the same time as Fefferman's work on the ball multiplier, the two-dimensional version of the
Bochner-Riesz conjecture was fully resolved
by Carleson and Sj\"olin~\cite{CarlesonSjolin}. 
More than a half-century later, though, the
 higher dimensional version of the conjecture remains open,
highlighting the importance of ensuring that Kakeya sets in $\Rn$, $n\ge 3$, cannot be ``too small''.

\subsection{The Kakeya Conjecture}
Given Fefferman's 1971 results, there has naturally been much interest since then in the following
modern form of the Kakeya problem:

\begin{conjecture}[Kakeya Conjecture]  Let $n\ge 2$.  Then if ${\mathcal K}$ is a Kakeya set in $\Rn$,
\begin{equation}\label{2.14}
\dim {\mathcal K}=n,
\end{equation}
with $\dim$ denoting either Minkowski dimension, $\dim_M$, or, more importantly,
Hausdorff dimension, $\dim_H$.
\end{conjecture}

We refer the reader to \cite{Fal} for precise definitions of Minkowski and Hausdorff dimensions.  We recall, though,
that $\dim_H E\le \dim_M E$ for any $E\subset \Rn$, which means that the Hausdorff form of the Kakeya conjecture
is stronger.  To avoid technicalities, we shall mainly focus on the Minkowski version of the conjecture. 

One form of the Minkowski variant of \eqref{2.14} is easy to state.  Indeed, if 
\begin{equation}\label{2.15}
{\mathcal K}_\delta =\{\, y\in \Rn: \, \text{dist }(y,{\mathcal K})<\delta\, \},
\end{equation}
denotes a $\delta$-neighborhood of a Kakeya set ${\mathcal K}$, then \eqref{2.14} with $\dim =\dim_M$ is just the following lower
bound for its measure:
\begin{equation}\label{2.16}
\text{For every } \, \, \e>0, \, \, \, |{\mathcal K}_\delta|\ge c_{\e,{\mathcal K}} \, \delta^\e, \quad
\text{whenever } \, \, \delta \in (0,1/2),
\end{equation}
with $c_{\e,{\mathcal K}}>0$ depending on $\e$ and the Kakeya set.  Note that if ${\mathcal K}$ is a zero-measure Kakeya set,
then by the dominated convergence theorem, $|{\mathcal K}_\delta|\to 0$ as $\delta\to 0$.   The lower bound
\eqref{2.16} says that the measure cannot go to zero too quickly.

An equivalent formulation of the statement $\dim_M{\mathcal K}=n$
 is that if
$N(\delta)$ is the minimum number of $\delta$-balls needed to cover ${\mathcal K}_\delta$, then
$\lim_{\delta\to 0} \log N(\delta)/\log \delta^{-1}=n$.   
The definition of Hausdorff dimension involves countable coverings by
sets of diameter $\delta$ or less, and this extra flexibility accounts for why $\dim_H E\le \dim_M E$.

In 1971, the same year as Fefferman's important result, R. Davies~\cite{RoyD} completely resolved the
two-dimensional version of the conjecture.

\begin{theorem}[Davies 1971]  If ${\mathcal K}\subset \R^2$ is a Kakeya set then
$\dim {\mathcal K}=2$.
\end{theorem}

Since Davies’ work, many proofs of this two-dimensional result have been obtained.  The simplest are $L^2$-based.  Indeed,
by using ideas going back to C\'ordoba~\cite{CDuke} (also from the 1970s), one can prove that Kakeya subsets of $\R^2$ always have full dimension by a simple argument involving the Cauchy-Schwarz inequality.  Davies was a student of 
Besicovitch, and C\'ordoba was a student of Fefferman.

\section{Wang and Zahl's Resolution of the three-dimensional Kakeya Conjecture}\label{kaksec}

We now turn to the work of Hong Wang.  One of her main results is the third milestone in the Kakeya story.
In 2025,  she  and Joshua Zahl were able to fully resolve the three-dimensional Kakeya conjecture in the remarkable work \cite{WZ3}.
As we shall discuss below, this is the first of the principal conjectures in harmonic analysis related to 
the Fourier transform to be settled in a dimension larger than two.

The result of Wang and Zahl, coming more than a half century after its two-dimensional analog, is the following.

\begin{theorem}[Wang and Zahl, 2025]\label{thm3.1}
If ${\mathcal K}$ is a Kakeya subset of $\R^3$, then it must have full dimension, both Minkowski and
Hausdorff.  Thus, $\dim_M{\mathcal K}=\dim_H{\mathcal K}=3$.
\end{theorem}

We shall attempt to give some of the main ideas  in the proof of this remarkable result.  The proof was
split up over three papers, \cite{WZ2}, \cite{WZ3} and \cite{WZ1},  spanning more than 300 pages.

To minimize technicalities, we shall focus on the special case saying that Kakeya subsets of $\R^3$ always
have full Minkowski dimension.  As in the earlier discussion of background results, we shall
also take a great deal of poetic license in places, for instance ignoring the role of  various small
constants that necessarily arise in the proof.

Given the importance of the result, not surprisingly, there have been several expository articles about the work
of Wang and Zahl.  These include Guth~\cite{G1}, \cite{G2},  \cite{G3}, Hickman~\cite{Hick} and
Tao~\cite{TKak}.  There is also the excellent article and accompanying
{\color{red}\href{https://www.youtube.com/watch?v=5J3tYU_-IZI}{video}} 
  by Quanta magazine
\cite{QKak}.

Before saying a few words about the proof of the Kakeya theorem of Wang and Zahl, let us mention some
earlier work on the Kakeya conjecture for $n\ge3$.  Given its importance, it is not surprising that many of the
leading harmonic analysts of the last couple of generations had made important contributions.  Many of the tools
developed to obtain partial results played a role in Wang and Zahl's proof.

Let us briefly mention  partial results for the three-dimensional Kakeya conjecture that preceded
Theorem~\ref{thm3.1}.  First on the list is the important 1991 work of Bourgain~\cite{B91} showing that
$\dim_H{\mathcal K}\ge 7/3$.  This was the first result going beyond $L^2$-methods.  Shortly 
afterwards, Wolff~\cite{W25} improved the lower bound to $\dim_H{\mathcal K}\ge 5/2$, and his approach
highlighted the power of multiscale analysis.  The Wolff lower bound proved to be difficult to improve upon
in the ensuing years.  The first improvement was by Katz, {\L}aba and Tao~\cite{KLT}, who showed
that $\dim_M {\mathcal K}\ge 5/2 + \e_0$ for some unspecified $\e_0>0$ very 
small (think $\e_0=10^{-10}$).  Though this improvement over the Wolff lower bound was modest (and only
for Minkowski dimension), this paper introduced very important ideas, including the notion of sticky Kakeya
sets.  
%This was followed by a related influential blogpost \cite{Tblog}.  
More recently,
Katz and Zahl~\cite{KZ} obtained a result in the spirit of \cite{KLT}, showing that 
$\dim_H {\mathcal K}\ge 5/2+\e_0$, also for some small unspecified $\e_0>0$.

\subsection{Discretization: Estimating unions of tubes}\label{discrete}
As is standard, the first step in the proof of Theorem~\ref{thm3.1} is to reduce to a discrete problem
involving $\delta$-tubes, $\{\Td\}=\{T\}$, meaning that each $T$ is a tube in $\Rt$ of length
$1$ and diameter $\delta \in (0,1/2)$.  We shall say that $\Td$ has $\delta$-separated directions
if the directions of the tubes $T$, $\omega_T\in S^2$, satisfy $|\omega_T-\omega_{T'}|\ge \delta$ for
all $T\ne T'$ in $\Td$.

We can now state the main estimate of Wang and Zahl to handle Minkowski dimension:

\begin{theorem}[Discrete Kakeya Theorem]\label{thm3.2}
Given $\e>0$ there is a constant $c_\e>0$ so that whenever $\Td=\{T\}$ is a collection of
$\delta$-separated $\delta$-tubes with $\delta\in (0,1/2)$, we have
\begin{equation}\label{3.1}
\bigl| \, \bigcup_{T\in \Td} T\, \bigr| \ge c_\e \delta^\e \, \bigl( \, |\Td|\cdot |T|\, \bigr).
\end{equation}
\end{theorem}

Here $|\Td|$ denotes the cardinality of the collection of tubes and $|T|\approx \delta^2$ the volume of 
each tube in $\Td$.  Although the discrete Kakeya estimate \eqref{3.1} of Wang and Zahl is very simple
(and trivial if the tubes happen to be disjoint), it is very difficult to prove.  One must consider {\em all possible
configurations} of directionally separated tubes, and, as Figure~\ref{fig3} suggests, there can be many intersections
potentially reducing $|\bigcup T|$.  The point of \eqref{3.1} is that $\delta$-separated tubes 
``cannot intersect too much''.

\begin{figure}[htbp]
  \centering
  \includegraphics[width=0.5\columnwidth]{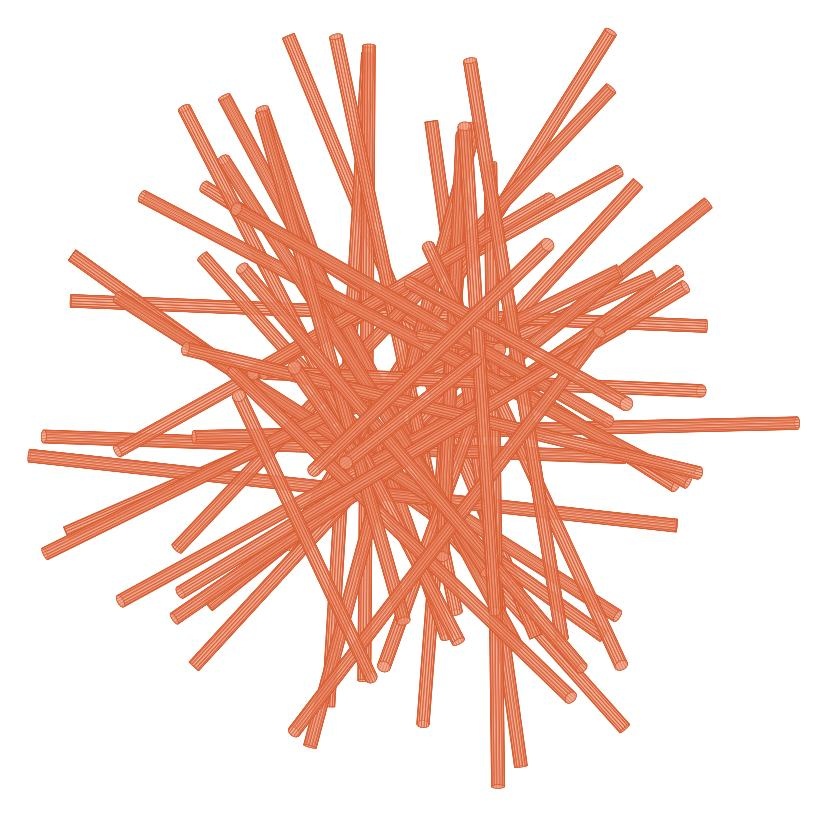}
  \caption{Union of $\delta$-tubes}
  \label{fig4}
\end{figure}

Let us give the simple argument showing that the discrete Kakeya estimates imply that $\dim_M {\mathcal K}=3$ if 
${\mathcal K}\subset \Rt$ is a Kakeya set.  We first note that, if as in \eqref{2.15}, ${\mathcal K}_\delta$ is a
$\delta$-neighborhood of ${\mathcal K}$, then, since, by assumption, ${\mathcal K}$ contains a unit line segment
in every direction, given $\omega\in S^2$ we can find a $\delta$-tube which points in the direction $\omega$
and belongs to ${\mathcal K}_\delta$.  Using the fact that $S^2$ is two-dimensional, we deduce from this that we
can find a collection, $\Td$, of $\delta$-separated $\delta$-tubes of cardinality $|\Td|\approx \delta^{-2}$ with
$T\subset {\mathcal K}_\delta$, $\forall T\in \Td$.  Whence, by \eqref{3.1},
$$|{\mathcal K}_\delta|\ge \bigl| \, \bigcup_{T\in \Td}T \, \bigr| \ge c_\e \delta^\e \,
\bigr(\, \delta^{-2}\cdot \delta^2\, \bigr)=c_\e \delta^\e, \, \, \, \forall \, \e>0.
$$
Thus, \eqref{2.16} must be valid and so our Kakeya set must satisfy $\dim_M {\mathcal K}=3$ as claimed.

To prove the stronger and more important result concerning Hausdorff dimension, $\dim_H {\mathcal K}=3$, Wang and Zahl
use a technical variant of \eqref{3.1}  involving ``shaded'' portions of each of the tubes.  Specifically, given
potentially small $\la>0$, they consider measurable subsets $Y(T)\subset T$ of each tube satisfying
\begin{equation}\label{3.2}
|Y(T)|\ge \la |T|.
\end{equation}
Given any such ``shading'' they are able to show that, for each $\e>0$ there is a $K\in {\mathbb N}$ so that,
for sufficiently small $\delta>0$
\begin{equation}\label{3.3}
\bigl| \, \bigcup_{T\in \Td} Y(T)\, \bigr| \ge \delta^\e \, \la^K\, 
\bigl( \, |\Td|\cdot |T|\, \bigr),
\end{equation}
assuming, as in \eqref{3.1}, that the tubes in $\Td$ are $\delta$-separated.

By a more technical argument than above, \eqref{3.3} implies that Kakeya subsets of $\Rt$ must have full
Hausdorff dimension.

\subsection{Sticky Kakeya sets: A significant special case}
Arguably the most significant step in Wang and Zahl's proof of the Kakeya conjecture was handled in
the first of their trilogy of papers, \cite{WZ1}.  Here they showed that the conclusions of Theorem~\ref{thm3.2}
are valid in the {\em sticky case}.

\begin{figure}[htbp]
  \centering
  \includegraphics[width=0.9\columnwidth]{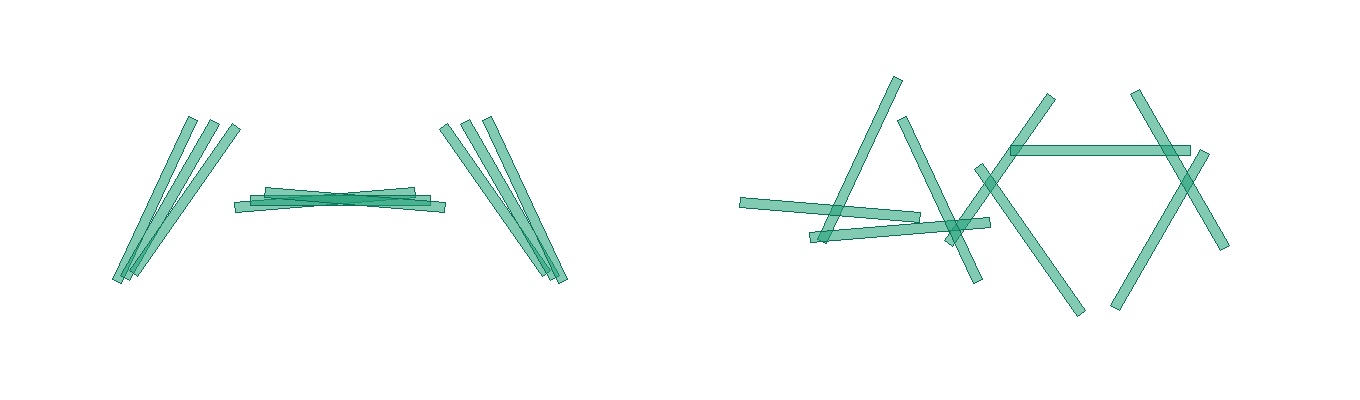}
  \caption{Sticky vs Non-Sticky in 2-d}
  \label{fig5}
\end{figure}

Informally, a collection of $\delta$-tubes, $\Td$, is sticky if typically tubes pointing in similar directions must also
sit together close in space.  So, as in Figure~\ref{fig5}, in the sticky case tubes with nearby directions cluster in
nearby locations (like branches of a tree), while for a non-sticky collections, tubes can share nearly the same direction
and sit in totally unrelated parts of space.

Due to the multiscale arguments used, it is convenient to formulate things more precisely using coverings of the
tubes $\Td$ by tubes $\Tr$ of coarser intermediate scales $\rho\in (\delta,1/2)$.  Specifically, we say that a collection
of $\delta$-separated $\delta$-tubes $\Td$ is sticky if, for every such $\rho$, we can find a collection
of $\rho$-tubes $\Tr$, with the following properties:
% i)  every $\delta$-tube $T_\delta\in \Td$ belongs to a
%$\rho$-tube $T_\rho$, ii) $|\Tr|=O(|\Td| \cdot (\delta/\rho)^2)$  , iii) if $\Td[T_\rho]=\{T\in \Td: \, T\subset T_\rho\}$, then
%$|\Td[T_\rho]|\approx (\rho/\delta)^2$ for a typical $T_\rho\in {\mathbb T}_\rho$.
\begin{enumerate}[label=\roman*.]
\item  every $\delta$-tube $T_\delta\in \Td$ belongs to a
$\rho$-tube $T_\rho$, 
\item $|\Tr|=O(|\Td| \cdot (\delta/\rho)^2)$,
\item and, if $\Td[T_\rho]=\{T\in \Td: \, T\subset T_\rho\}$, then
$|\Td[T_\rho]|\approx (\rho/\delta)^2$ for a typical $T_\rho\in {\mathbb T}_\rho$.
\end{enumerate}

\begin{figure}[htbp]
  \centering
  \includegraphics[width=0.9\columnwidth]{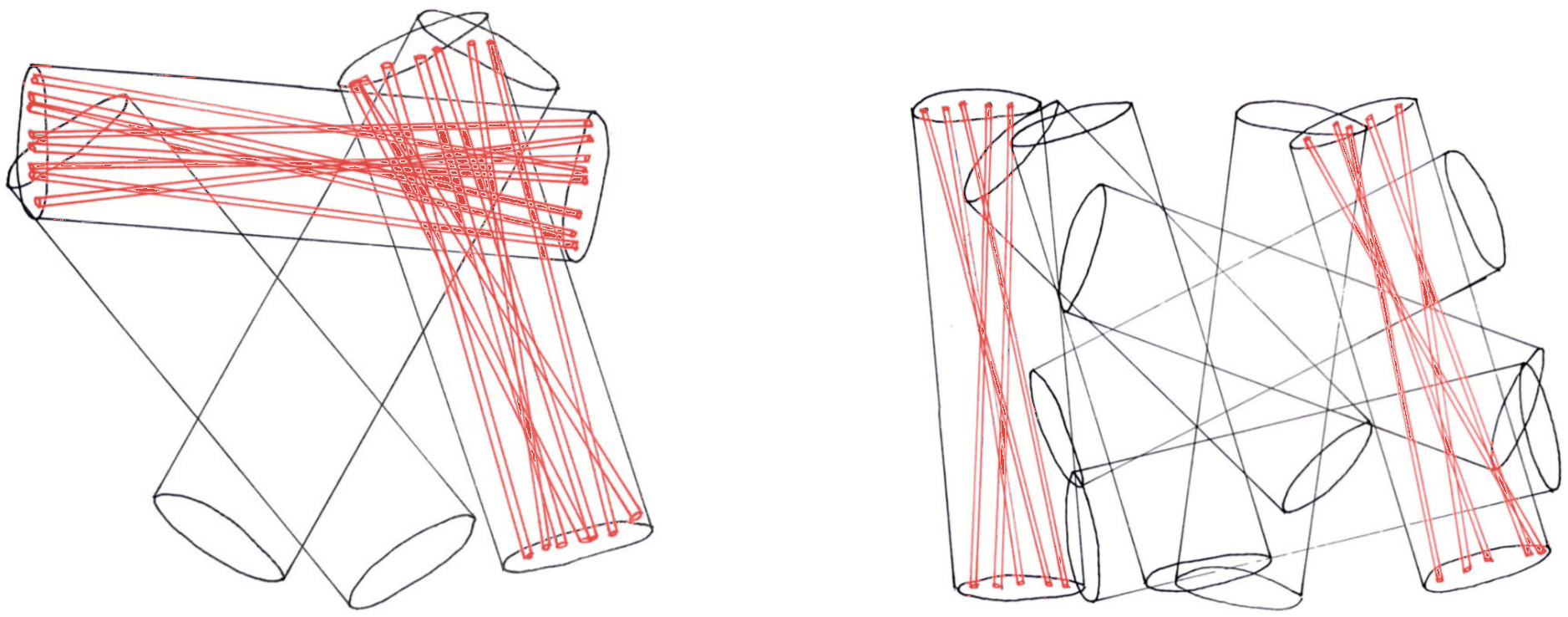}
  \caption{Sticky and Non-sticky collections of tubes in $\Rt$}
  \label{fig6}
\end{figure}

Since the set of directions in $\Rt$ is two-dimensional, one can think of a sticky collection of tubes, $\Td$, as
an urban planner's dream:  Given any scale $\rho\in (\delta,1/2)$ $\Td$ is covered by exactly the right number
$\approx |\Td| \cdot (\delta/\rho)^2$ neighborhoods $\{T_\rho\in \Tr\}$, with typical neighborhoods 
$\Td[T_\rho]=\{T\in \Td: \, T\subset T_\rho\}$ having the natural size $|\Td[T_\rho]| \approx (\rho/\delta)^2$.
  Little to no flyover country.  Sticky and non-sticky collections of tubes are depicted in Figure~\ref{fig6}.  For visualization
  reasons, we have only drawn the smaller $\delta$-tubes (shaded red) in some of the larger $\rho$-tubes.
  
  We shall say that a Kakeya set ${\mathcal K}\subset \Rt$ is {\em sticky} if for every $\delta\in (0,1/2)$,
  ${\mathcal K}_\delta$ (as above) can be covered by a collection $\Td$ of sticky $\delta$-tubes.  Then the main
  result from Wang and Zahl~\cite{WZ1} is the following.
  
  \begin{theorem}[Sticky Kakeya Theorem]\label{thm3.3}  Sticky Kakeya subsets of $\Rt$ have full dimension.
  Additionally, \eqref{3.1} is valid if the $\delta$-direction separated collection $\Td$ is sticky.
  \end{theorem}
  
  We should point out that the sticky condition is rigid.  Sticky Kakeya sets have a self-similar like structure.  For instance, if 
  one considers the $\delta$-tubes lying in a typical $\rho$-tube as in Figure~\ref{fig6} and performs an anisotropic scaling sending
  the $\rho$-tube to a 1-tube, the  image of the $\delta$-tubes become a collection of
$\sim (\delta/\rho)$-tubes that is also sticky.
  Wang and Zahl exploited the rich structure of sticky Kakeya
  sets in their proof.  On the other hand, they had to overcome two serious obstacles.  The first is intuitive.  
  Sticky Kakeya sets should be ones with the maximal amount of ``compression'' and naively the worst case
  for Theorem~\ref{thm3.1}.
  
  A more serious obstacle was pointed out by {\L}aba, Katz and Tao~\cite{KLT}: the analog of Theorem~\ref{3.2}
  is false in the {\em complex setting!}  
%  Hence, any proof of the Sticky Kakeya Theorem must distinguish between
%  the reals and ${\mathbb C}$.
 
 Briefly, the ``Heisenberg example'' from \cite{KLT} is given by
  \begin{equation}\label{3.4}
  {\mathbb H}=\bigl\{(z_1,z_2,z_3): \, \, \text{Im }(z_3)=\text{Im }(z_1\, \overline{z}_2), \, \, |z_1|, \, |z_2|, \, |z_3|\le 2 \, \bigr\} \subset  {\mathbb C}^3.
  \end{equation}
  Note that ${\mathbb H}$ contains a 4 (real)-parameter family of line segments, $\{(z,w+az, z\overline{w}+b): \, |z|\le 1/2\}$, $a,b\in \R$, $w\in {\mathbb C}$ with
  $|a|, |b|, |w| \le 1$.  From this, one can see that there are $\sim \delta^{-4}$ complex $\delta$-tubes about complex unit segments, ${\mathbb T}$, with
  $\delta$-separated directions which
  belong
  to a $\delta$-neighborhood, ${\mathbb H}_\delta$, of ${\mathbb H}$.  Thus, $|\bigcup_{T\in {\mathbb T}}T| \lesssim \delta$.
  Since each complex tube has measure $\sim \delta^4$, it follows that $ |\bigcup T| \lesssim \delta  \ll 1\approx (|{\mathbb T}|\cdot |T| )$, which means the analog
  of \eqref{3.1} cannot hold here.
  
 The Heisenberg example also shares the key structural properties of sticky Kakeya subsets.
  From this we conclude that any proof of Theorem~\ref{thm3.3} must distinguish between the real and complex field.  The fundamental
  structure of the Heisenberg example relies on the fact that ${\mathbb C}$ has the half-dimensional subfield ${\mathbb R}\subset {\mathbb C}$.  
  It turns out that ${\mathbb R}$ cannot have subfields of intermediate dimension.  Indeed, an important theorem of
  Edgar and Miller~\cite{EM} says that Borel subrings of ${\mathbb R}$ must either be zero-dimensional or be all of ${\mathbb R}$.
  
  In order to prove their sticky Kakeya theorem Wang and Zahl used results that grew out of the celebrated sum-product theorem
  of Bourgain~\cite{Bsp}, which can be thought of as a quantitative variant of Edgar and Miller's theorem.  More precisely, Wang and Zahl 
  used  recent related results from projection theory (\cite{32}, \cite{34}, \cite{35}) that were developed in the study of the Falconer distance problem.
In  particular, the results of Orponen, Shmerkin and Wang~\cite{35}
  help to provide the critical step in the proof which distinguishes between a sticky Kakeya set in $\Rt$  and the Heisenberg
  example (a subset of ${\mathbb C}^3$).
  
  \subsection{End of proof: Leveraging sticky Kakeya}
  
  The directional hypotheses in the Discrete Kakeya Theorem (Theorem~\ref{thm3.2}) are not typically preserved in the multiscale arguments
  that Wang and Zahl employ.  Indeed, after zooming in on a small region of $\bigcup_{T\in \Td} T$ and rescaling, one could end up with a new collection
  of larger tubes 
  %having  many 
  where many of the tubes are
  pointing in nearly the same direction.  

Due to this, Wang and Zahl were led to proving stronger estimates involving hypotheses that can be used in their multiscale induction arguments.
The assumptions take into account the density of $\Td$ in convex subsets ${\mathcal C}$ in $\Rt$.  Specifically,
\begin{equation}\label{235}
\Delta(\Td, \mC)=\frac{\# \{T\in \Td: \, T\subset \mC\} \times \delta^2}{|\mC|} \approx
\frac{\text{combined volume of tubes in } \, \mC}{|\mC|}.
\end{equation}

\begin{figure}[htbp]
  \centering
  \includegraphics[width=0.8\columnwidth]{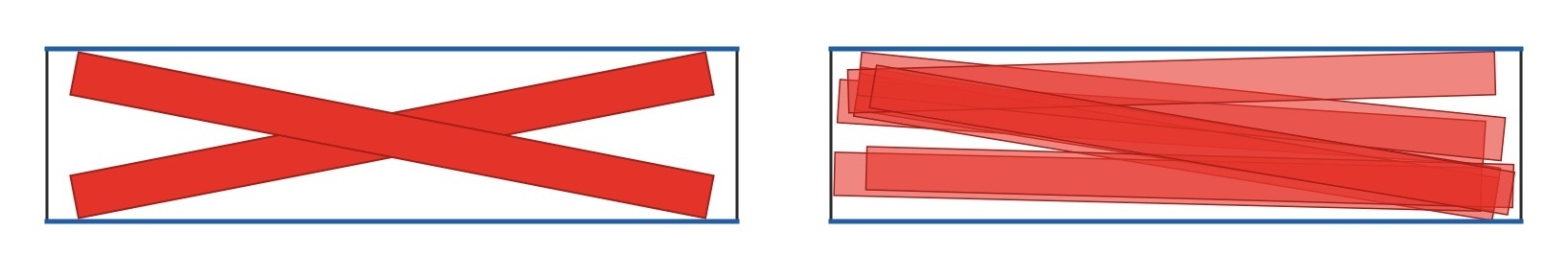}
  \caption{Low and high densities}
  \label{fig7}
\end{figure}

Due to the John ellipsoid theorem, one can always assume that the sets $\mC$ are rectangular prisms for the sake of simplicity, but this is not 
very important.  Figure~\ref{fig7} depicts low and high density clustering.

For the remainder of the section, we shall assume that all of the tubes $T\in \Td$ that arise are contained in the unit ball
$$B=\{x\in \Rt: \, \, |x|\le 1\}.$$
The stronger versions of the earlier discrete Kakeya theorem are based on two different types of tube clustering conditions:

\begin{definition}\label{KT}
A collection of $\delta$-tubes $\Td$ satisfies the {\bf Katz-Tao} condition if
\begin{equation}\label{3.6}
\sup_{\mC\subset \Rt, \, \, \text{convex}} \Delta(\Td, \mC)\lesssim 1.
\end{equation}
\end{definition}

\begin{definition}\label{F}
A collection of $\delta$-tubes $\Td$ satisfies the {\bf Frostman} condition if
\begin{equation}\label{3.7}
\sup_{\mC\subset \Rt, \, \, \text{convex}} \Delta(\Td, \mC)\lesssim 
\Delta(\Td,B).
\end{equation}
\end{definition}

Wolff seems to be one of the first to have organized tubes as above.  In his important Kakeya paper \cite{W25}, his main result was formulated using the  assumption \eqref{3.6}, and in
\cite{Wx} he introduced the other type of clustering condition \eqref{3.7}.  For these reasons, Wang and Zahl also call \eqref{3.6} the {\em Katz-Tao Convex Wolff Axiom}
and \eqref{3.7} the {\em Frostman Convex Wolff Axiom}.  Wolff's paper \cite{Wx} concerned the X-ray transform, while the papers
of Wang and Zahl seem to be the first to use the Frostman condition for the Kakeya problem.

It is important to  emphasize that the two conditions are complementary.  For instance, the Katz-Tao condition is an absolute condition, imposing an absolute
bound on the clustering density, while the Frostman condition is a relative condition saying that the clustering of tubes within convex sets is dominated
by the clustering within the ball, $B$.  The two conditions also give complementary requirements for the cardinality of $\Td$.  Indeed, by taking
$\mC=B$ in \eqref{3.6} we conclude that if the Katz-Tao condition holds then $|\Td|\lesssim \delta^{-2}$, while, conversely, by taking $\mC$ to be
a single tube $T$, we deduce that  if the Frostman condition holds then $|\Td|\gtrsim \delta^{-2}$.  Another important observation is that
if, as in Theorem~\ref{thm3.2}, the tubes $\Td$ are $\delta$-direction separated then  the Katz-Tao condition holds,
and the Frostman condition is  valid if additionally $|\Td|\gtrsim \delta^{-2}$.

It turns out that there is a very surprising 
{\em ``Yin and Yang''} relationship between the two conditions that Wang and Zahl are able to exploit in their remarkable inductive proof of the following results.

\begin{theorem}(Strong Kakeya: Katz-Tao version)\label{thmkt}
If $\Td$ as above is {\em Katz-Tao}, then \eqref{3.1} is valid,   i.e.,
%i.e., there is a constant $c_\e>0$ so that
%for every $\e>0$,
\begin{equation}\label{3.8}
\bigl| \, \bigcup_{T\in \Td} T\, \bigr| \ge c_\e \delta^\e \, \bigl( \, |\Td|\cdot |T|\, \bigr), \quad \forall \, \e>0.
\end{equation}
\end{theorem}

\begin{theorem}(Strong Kakeya: Frostman version)\label{thmf}
If $\Td$ as above is {\em Frostman}, then 
%there is a constant $c_\e>0$ so that for every $\e>0$
\begin{equation}\label{3.9}
\bigl| \, \bigcup_{T\in \Td} T\, \bigr| \ge c_\e \delta^\e, \quad \forall \, \e>0.
\end{equation}
\end{theorem}

Clearly, the arguments of \S\ref{discrete} imply that either of these two theorems implies that Kakeya subsets of $\Rt$ must have
full dimension.  Also, of course, given the remarks before the statement, Theorem~\ref{thmkt} yields Theorem~\ref{thm3.1}, and 
Theorem~\ref{thmf} yields the important special case where $|\Td|\approx \delta^{-2}$.

The two theorems are proved using a coupled inductive argument which takes advantage of the complementary
nature of the two clustering conditions.  The two inductive estimates\footnote{We are presenting here
a slight variation from \cite{GWZst} of inductive bounds in \cite{WZ3}.} say that for 
$\beta\in [0,1]$
\begin{equation*}
\text{KT}(\beta):   \qquad \qquad 
\bigl| \, \bigcup_{T\in \Td} T \, \bigr|
\gtrsim \delta^\e \, |\Td|^{-\beta} \, \bigl( \, | \Td| \cdot |T|\, \bigr), \, \, \forall \, \e>0, 
\qquad \text{if } \, \Td \, \, \text{is Katz-Tao},
\end{equation*}
and, additonally,
\begin{equation*}
\text{F}(\beta):   \qquad \qquad 
\bigl| \, \bigcup_{T\in \Td} T \, \bigr|
\gtrsim \delta^\e \, \delta^{2\beta} \,  \bigl( \, | \Td| \cdot |T|\, \bigr)^{\beta/2}, \, \, \forall \, \e>0, 
\qquad \text{if } \, \Td \, \, \text{is Frostman}.
\end{equation*}
We are oversimplifying the bounds a bit here for reasons of exposition.  Wang and Zahl prove (and their
arguments require) stronger estimates which are a bit more technical in nature, involving, for instance
the shading described in \S3.1.

Note that since $|T|\approx \delta^2$, it is easy to see that both KT($\beta$) and F($\beta$) are valid when
$\beta=1$.  Also, by the arguments in \S3.1, for the special case where $|\Td|\approx \delta^{-2}$ and the tubes
have $\delta$-separated directions, either of the two estimates imply that Kakeya subsets of $\Rt$ must have
dimension $\ge 3-2\beta$.  Thus KT($\beta$) and F($\beta$) can be thought of as partial Kakeya estimates.

Here comes the striking argument of Wang and Zahl.  Oversimplifying a bit (e.g., ignoring shading issues), they
are able to carry out a self-improving inductive argument showing that: 
%KT($\beta$) implies F($\beta$) {\em and} if both KT($\beta$) and F($\beta$) are valid then 
%KT($\beta-\nu(\beta)$) must also be valid for some $\nu(\beta)>0$.  
\begin{enumerate}[label=\roman*.]
\item KT($\beta$) implies F($\beta$), {\em and} 
\item  if  $\beta>0$ and both KT($\beta$) and F($\beta$) are valid,  then 
KT($\beta-\nu(\beta)$) must also be valid for some $\nu(\beta)>0$. 
\end{enumerate}
Since the estimates include the factor
$\delta^\e$ in the right, this means that the set of $\beta\in [0,1]$ for which the two estimates hold
must be open and closed in $[0,1]$.  This set is nonempty because, as mentioned before, 
KT(1) and F(1) are both valid.  Thus, KT(0) and F(0) are also valid, and hence Theorems~\ref{thmkt}
and \ref{thmf} must  also hold!!

To summarize, Wang and Zahl set up an inductive argument that involves two complementary 
variations on the statement that Kakeya sets ${\mathcal K}\subset \Rt$  must satisfy
$\dim {\mathcal K}\ge 3-2\beta$.  Of course this is true for $\beta=1$ because trivially
$\dim {\mathcal K}\ge 1$ as the Kakeya set contains line segments.  Then, through multiscale arguments
that we shall briefly discuss, they use the Yin and Yang nature of the two estimates to show  that
the dimension has to be a bit larger until concluding that Kakeya sets always of dimension larger than
$3-\delta$ for any $\delta>0$
after finitely many steps!

We have attempted to describe, in very rough terms, the remarkable inductive argument used
by Wang and Zahl to resolve the three-dimensional Kakeya conjecture.  Let us also try to sketch
how they are able to leverage their sticky Kakeya theorem to run the inductive argument just
described.  They developed the necessary tools in the last two papers \cite{WZ2} and \cite{WZ3}
of their trilogy.  Here they strengthened their ``sticky'' results and also found novel ways of organizing
tubes to take advantage of structures occurring in the various steps of their multiscale arguments.

To carry out the above inductive argument, Wang and Zahl take advantage of the fact that, if
$\beta>0$, they are able to prove stronger estimates than KT($\beta$) and F($\beta$) 
in the sticky case.  Therefore, they can reduce to a range of scales for which the sticky criteria
strongly break down.  This can be expressed using Katz-Tao or Frostman clustering bounds.
For the former, it means that if ${\mathbb T}_\rho$ is a collection of $\rho$-tubes covering the
$\delta$-tubes, then, when $\delta$ is replaced by $\rho$ in \eqref{3.6}, the resulting sup is very
large.  This case, as in the right half of Figure~\ref{fig6}, is an urban planner's
nightmare, with a large number of $\rho$-tubes needed to cover $\Td$.

\begin{figure}[htbp]
  \centering
  \includegraphics[width=0.6\columnwidth]{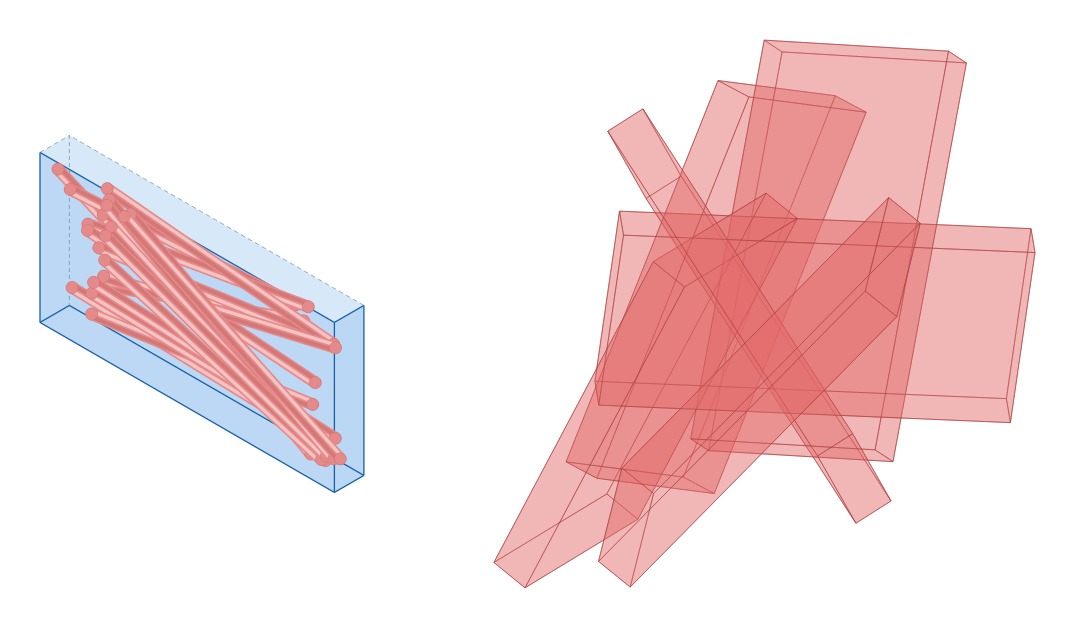}
  \caption{$\Td[R]$ and  ${\mathcal R}=\{R\}$ }
  \label{fig8}
\end{figure}

To take advantage of this, Wang and Zahl use their organizational skills.  They establish
a ``factorization lemma.''  Using a straightforward pigeonhole argument they show that,
after possibly culling the collection of tubes $\Td$ a negligible amount, they can find a collection
${\mathcal R}=\{R\}$ of $a \times b\times 1$ rectangles, with $a\le b\le 1$ so that $\Td$ can be written,
as depicted in Figure~\ref{fig8},
as the disjoint union
\begin{equation}\label{33.1}
\Td = \bigsqcup_{R\in {\mathcal R}} \Td[R],
\quad \text{with } \, \, \, \Td[R]=\{T\in \Td: \, \, T\subset R\}.
\end{equation}
where:
% i) the tubes $T\in \Td[R]$ are Frostman in $R$, and
%ii) the $a\times b\times 1$ rectangles ${\mathcal R}$ are Katz-Tao
%in $B$.  
\begin{enumerate}[label=\roman*.]
\item  the tubes %$T\in \Td[R]$ 
$\Td[R]$ 
are Frostman in $R$, and
\item the $a\times b\times 1$ rectangles ${\mathcal R}$ are Katz-Tao
in $B$.  
\end{enumerate}
The first condition means that the analog of \eqref{3.7} is valid with $B$
replaced by $R$, and ii) means that the analog of \eqref{3.6} is  valid
where $\Td$ is replaced by ${\mathcal R}$ (and in \eqref{235} the densities are
similarly modified).

Due to the factorization \eqref{33.1}, the estimates for
$|\cup T|$ in KT($\beta$) and F($\beta$) can be obtained by obtaining appropriate lower bounds
for the union of the tubes in each $a\times b\times 1$ rectangle $R$, i.e., $|\cup_{T\in \Td[R]}T|$,
using the Frostman assumption, and then controlling volume lower bounds when unions of these
sets are taken over ${\mathcal R}$ using the Katz-Tao assumption for the rectangles.

Thus, Wang and Zahl reduce to controlling the interactions of Frostman tubes within $a\times b\times 1$
rectangles, and the interactions of shaded Katz-Tao rectangles with these side lengths.  
They are able to use the inductive bounds  KT($\beta$) and F($\beta$)  together with a simple organizing of rectangles according to planar direction argument to prove natural generalizations of these two estimates where
tubes are replaced by
rectangles for the variant of KT($\beta$), and all the tubes are assumed to lie in an $a\times b\times 1$ rectangle for
the F($\beta$) variant.

%The {\em tour de force} argument yielding these bounds would lead to nothing without assumptions on the shape
%of the rectangles ${\mathcal R}$.  We can finally say how they leverage their sticky results.  

Recall that to carry out a difficult step of
the self-improving inductive argument they can assume that $\Td$ is {\em very non-sticky}.  Miraculously,
after dealing with some relatively straightforward cases, they are able to use this to their advantage by showing that each
$R\in {\mathcal R}$ in the key step of the proof must be shaped like an iPhone, by which we mean that $a\ll b<1$.  So, the 
``short side'' of $R$ must be much smaller than the ``medium side'' due to the very non-sticky assumption.
The variants of KT($\beta$) and F($\beta$)  involving ${\mathcal R}$ that they obtain
in this part of the inductive argument
 involve positive powers of $b/a\gg 1$ in the right, and, as a result,
they can cash in by using the auxiliary bounds involving ${\mathcal R}$ to carry out their self-improving 
%inductive
argument!

\section{Principal conjectures related to the Fourier transform}

Let us discuss Hong Wang's work on other problems related to the Fourier transform.

Earlier we presented the Bochner-Riesz and Kakeya conjectures.  There are two other principal conjectures with close
connections to the Fourier transform.  These four conjectures form a family of problems with surprising connections.
As we shall see, Hong Wang has made very important contributions to all four.

Of the two not yet discussed, the oldest is Stein's restriction problem dating back to the 1970s:
\begin{equation}\label{4.1}
\text{Show that for each } \, n\ge 2,  \quad
\| \Hat f \|_{L^q(S^{n-1})} \le C_p \|f\|_{L^p(\Rn)}, \, \, f\in {\mathcal S}(\Rn), \, \, \,
\text{if } \, \, 1\le p<\tfrac{2n}{n+1} \, \, \text{and } \, \, q=\tfrac{n-1}{n+1} p'.
\end{equation}
Zygmund~\cite{Zyg} in 1974 fully resolved this problem for the two-dimensional case.  Thus, three of the four conjectures 
were fully understood in $\R^2$ in the early 1970s.

The other main conjecture was formulated by the author \cite{Sls} in 1991.  It is a conjecture about optimal local $L^p$
space-time regularity for the wave equation,
\begin{equation}\label{4.2}
\begin{cases}
\bigl(\parital_t^2-\Delta\bigr) u(x,t)=0
\\
u(x,0)=u_0(x), \quad \partial_t  u(x,0)=u_1(x),
\end{cases}
\end{equation}
with $n\ge2$.  To state the conjecture, let $L^p_\delta(\Rn)$ be the $L^p(\Rn)$-Sobolev space, associated with $\delta$-derivatives, with norm
$$\|f\|_{L^p_\delta(\Rn)}= \bigl\| \, (I+P)^{\delta} f\|_{L^p(\Rn)}, \quad P=\sqrt{-\Delta}.$$
Then the local smoothing conjecture says that for $u$ as in \eqref{4.2} we have
\begin{equation}\label{4.3}
\|u\|_{L^p(\Rn\times [1,2])} \lesssim_\e \, \|u_0\|_{L^p_{\delta(p,n)+\e}(\Rn)}+\| u_1\|_{L^p_{\delta(p,n)-1+\e}(\Rn)}, \, \forall \, \e>0,
\quad \text{if } \, \, p\ge \tfrac{2n}{n-1} \, \, \,  \text{and } \, \, \delta(p,n)=n(\tfrac12-\tfrac1p)-\tfrac12, 
\end{equation}
or equivalently, for $P$ as above, and all $\e>0$,
\begin{equation}\label{4.4}
\|e^{itP}f\|_{L^p(\Rn\times [1,2])} \lesssim_\e \, \|f\|_{L^p_{\delta(p,n)+\e}(\Rn)}, \quad
\text{with } \, p \, \, \text{and } \, \, \delta(p,n) \, \, 
\text{as in } \, \, \eqref{4.3}.
\end{equation}
It is straightforward to obtain the optimal fixed-time estimates
$e^{itP}: L^p_{(n-1)(1/2-1/p)}(\Rn)\to L^p(\Rn)$, if $t\in [1,2]$ and $p\ge 2$ (see \cite{Peral}).
The estimates \eqref{4.3} and \eqref{4.4} say that for $p\ge \tfrac{2n}{n-1}$ one essentially gains $1/p$ derivatives by taking local
space-time norms, and it is also easy to see that this  would be the optimal range of exponents with this property.  Note that
the index $\delta(p,n)$ in the above conjecture agrees with the one in the Bochner-Riesz conjecture \eqref{2.10}.

The $n$-dimensional local smoothing conjecture is difficult and interesting because it implies the other principal conjectures
for the same dimension.

In \cite{Sls} the simple argument showing that the local smoothing conjecture for the wave equation implies the Bochner-Riesz
conjecture was given.  This just boils down to the fact that the half-wave operators, $e^{itP}$, and the Bochner-Riesz operators
\eqref{BR} are related via the Fourier transform:
$$S^\delta =\int_{-\infty}^\infty e^{-it P} e^{it} m_\delta(t) \, dt, \quad \text{with } \,
m_\delta(t)=O\bigl(\, (1+|t|)^{-1-\delta} \, \bigr).
$$
Using this formula and H\"older's inequality, it is not too difficult to show that the local smoothing conjecture \eqref{4.4} implies
\eqref{2.11} for $p\ge \tfrac{2n}{n-1}$, which, in turn, implies the Bochner-Riesz conjecture for all exponents $p$ by duality
and a simple interpolation argument.

It turns out that all four conjectures are related.  In fact they form a hierarchy.  Tao~\cite{TaoBR} showed that the 
Bochner-Riesz conjecture implies the restriction conjecture.  Bourgain in \cite{B91} showed that the restriction conjecture implies
the Hausdorff version of the Kakeya conjecture.  Thus, the four conjectures form a ``tower'' for each dimension $n\ge2$:
\begin{equation}\label{4.5}
\text{Local smoothing conjecture} \implies 
\text{Bochner-Riesz conjecture } \implies 
\text{Restriction conjecture } \implies 
\text{Kakeya conjecture}.
\end{equation}

\subsection{Resolution of the planar local smoothing conjecture for the wave equation}

There was much work on the local smoothing conjecture for the wave equation.  Early work in
\cite{Sls} and \cite{MSSls} showed that there is some gain in regularity versus the fixed-time
estimates when taking space-time norms.  These results yielded new proofs of Bourgain's circular
maximal theorem \cite{Bcirc}, etc.; however, they were far from optimal.

A very important breakthrough came from Wolff~\cite{Wls}.  In order to obtain optimal local smoothing
estimates for very large $p$ when $n=2$, he introduced decoupling theory, which has had a wide
range of applications beyond the local smoothing conjecture \eqref{4.3}, including ones
in number theory \cite{BDG}.  
%Posthumously, after Wolff's untimely death, {\L}aba and Wolff~\cite{LWls}
%proved similar results in higher dimensions, and, in particular, establishing \eqref{4.3} for large $p$. 
The paper of  {\L}aba and Wolff~\cite{LWls}, completed and published posthumously after 
Wolff's untimely death,
extended the partial results in \cite{Wls} to higher dimensions.
%, and, in particular, establishing \eqref{4.3} for large $p$. 

Decoupling theory was developed by several authors in the years after Wolff's paper, culminating in
the optimal  $\ell^2$ decoupling bounds of Bourgain and Demeter~\cite{BDls}.

To describe the special case of their results which is related to the problem \eqref{4.3},  we need to consider
a portion of the light cone in $\Rn\times \R$:
$$\Gamma = \, \bigl\{ \, 
(\xi,|\xi|): \, \, |\xi|\in [1,2], \, \, \text{and } \, \, |\xi'|/|\xi|\le 1/2\, \bigr\}, \, \, \, \text{with } \, \, \xi'=(\xi_1, \dots,\xi_{n-1}),
$$
and its $\delta$-neighborhood ($0<\delta \ll 1)$
$$\Gamma_\delta = \{ \, (\eta,\tau) \in \Rn\times \R: \, \, \text{dist}((\eta,\tau),\Gamma)<\delta, \, \, \tau\in [1,2] \}.$$
Finally, we let $\text{Part}_{\delta^{1/2}} [-1,1]^{n-1}$ denote a partition of the cube $[-1,1]^{n-1}\subset \R^{n-1}$ into
non-overlapping cubes $\{J\}$ of side length $\delta^{1/2}$.  To each, we associate projection operators
$$P_J F(x,t)=\int_{1}^2 \int_{\{\xi\in \Rn: \, \xi'/\tau \in J\}}
e^{2\pi i \langle \, (x,t), (\xi,\tau) \, \rangle}
\Hat F(\xi,\tau) d\xi d\tau.$$
The decoupling estimates for $\Gamma$ of Bourgain and Demeter then say that for $n\ge2$
\begin{equation}\label{4.6}
\|F\|_{L^p(\Rn\times \R)}\lesssim_\e
\delta^{-\alpha(p)-\e}
\, \bigl( \, \sum_{J\in \text{Part}_{\delta^{1/2}} [-1,1]^{n-1}}
\|P_J F\|^2_{L^p(\Rn \times \R)} \, \bigr)^{1/2}, \quad
\text{if } \, \, \text{supp }\Hat F\subset \Gamma_\delta \, \, 
\text{and } \, \, p>2,
\end{equation}
with 
$$\alpha(p)=
\begin{cases}
\tfrac{n-1}4 (1-\tfrac{2(n+1)}{(n-1)p}) \quad \text{if } \, \, p\ge \tfrac{2(n+1)}{n-1}
\\
0\quad \text{if } \, \, 2<p\le \tfrac{2(n+1)}{n-1}.
\end{cases}
$$

\begin{figure}[htbp]
  \centering
  \includegraphics[width=0.25\columnwidth]{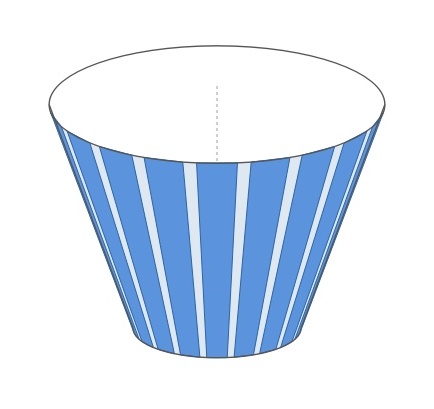}
  \caption{Cone decoupling }
  \label{fig9}
\end{figure}

The Fourier decomposition here is depicted in Figure~\ref{fig9}.  Each $P_JF$ is a ``wave packet'' with Fourier support
belonging to one of the shaded ``planks''.

An immediate corollary of the important decoupling bounds \eqref{4.6} is that \eqref{4.3} and \eqref{4.4} are valid for every
$p\ge \tfrac{2(n+1)}{n-1}$.  
They do not, however, lead to optimal bounds for smaller exponents.  For instance, when
$n=2$ and $p=4$ is the critical exponent, they merely imply the weaker local smoothing estimates from \cite{MSSls}
saying that there is local smoothing for $L^4$ of order $1/8-\e$ instead of the conjectured gain of
$1/4-\e$ derivatives.

Using novel multiscale ``wave envelope'' arguments Guth, Wang and Zhang \cite{GWZls} were able to prove sharp
square function estimates for $n=2$ that finally resolved the planar local smoothing conjecture.  Their square function bounds
superficially resemble \eqref{4.6}; however, it is much more difficult (and powerful) when $p=4$ and $n=2$ since the
$\ell^2$ norm is taken inside the $L^4$-norm:

\begin{theorem}[Local Smoothing Theorem of Guth, Wang and Zhang]\label{thmgwz}   If $n=2$, then
\begin{equation}\label{4.7}
\|F\|_{L^4(\R^2\times \R)} \lesssim_\e   \delta^{-\e} 
\bigl\| \, \bigl(\sum_{J\in \text{Part}_{\delta^{1/2}} [-1,1]} |P_JF|^2 \, \bigr)^{1/2} \, \bigr\|_{L^4(\R^2\times \R)}
 \, , \, \forall \, \e>0, \, \, \,
\text{if } \, \, \text{supp } \Hat F\subset \Gamma_\delta.
\end{equation}
Consequently, both \eqref{4.3} and \eqref{4.4} are valid when $n=2$.
\end{theorem}

It was well known (see \cite{MSSls}) that the sharp square function estimate \eqref{4.7} implies the conjecture for
$p=4$ and $n=2$, and the conjecture for all of the other exponents follows from this and a simple interpolation argument.

We have explained in some detail Wang's work on the first and last of the conjectures listed in \eqref{4.5}.  She also has made
significant contributions to the restriction problem,
obtaining the best known results, starting with her 2019 thesis work under  Guth \cite{Wrest1}, \cite{Wrest2} and
continued with her work with  Wu \cite{Wrest3}.  Similarly, her work with Guo, Oh, Wu and Zhang \cite{WangBR} provided the
best known results for the Bochner-Riesz conjecture.

Let us conclude this section by emphasizing how in \cite{GWZls} Hong Wang and collaborators finally solved the only one
of the four principal conjectures in \eqref{4.5} related to the Fourier transform that was not fully understood in the early 
1970s.  Furthermore, in her joint paper \cite{WZ3}, she was able to resolve, for the first time, one of these conjectures in 
higher dimensions.  So, she closed the door in two dimensions and opened it in higher dimensions, which is quite
an accomplishment!

\section{Other works}  We have only given a small sampling of Hong Wang's work.
Besides her work on the Fourier transform, she has several other very important
results.
These include her joint work on the Falconer distance problem with Larry Guth, Alex Iosevich, and Yumeng Ou
\cite{F3}, and her
resolution of the Furstenberg Set Conjecture with Kevin Ren~\cite{RenW}.

\section*{Acknowledgments.}
The figures were all prepared using Anthropic's Claude AI.  
Additionally, Figures  \ref{fig3}, \ref{fig6} and \ref{fig9} were adapted
from ones in \cite{Hick}, \cite{WZ3} and \cite{GWZls}, respectively.

\bibliography{refs}
\bibliographystyle{abbrv}

\end{document}